\documentclass[11pt]{amsart}
\usepackage{amscd,amssymb}

\newtheorem{Thm}{Theorem}[section]

\newtheorem{Cor}[Thm]{Corollary}
\newtheorem{Lem}[Thm]{Lemma}
\newtheorem{Prop}[Thm]{Proposition}

\begin{document}

\title[Character Varieties and Harmonic Maps]%
      {Character Varieties and Harmonic Maps to ${\mathbb R}$-trees}

\author[Daskalopoulos]{G.  Daskalopoulos}

\address{Department of Mathematics \\
		Brown University \\
		Providence,  RI  02912}

\thanks{G.D. supported in part by NSF grant DMS-9803606 and the Institute
for Advanced
Study, Princeton, N.J.}

\email{daskal@gauss.math.brown.edu}
\author[Dostoglou]{S. Dostoglou}

\address{ Department of Mathematics \\ University of Missouri
\\Columbia, MO  65211}

\thanks{S.D. supported in part by the Research Board of the University of
Missouri and the Arts \& Science Travel Fund of the University of Missouri,
Columbia}

\email{stamatis@euler.math.missouri.edu}

\author[Wentworth]{R.  Wentworth}

\address{Department of Mathematics \\
   University of California \\
   Irvine,  CA  92697}

\email{raw@math.uci.edu}

\thanks{R.W. supported in part by NSF grant DMS-9503635 and a Sloan Fellowship}

\date{}

\begin{abstract}
We show that the Korevaar-Schoen limit of the sequence of equivariant
harmonic maps
corresponding to a sequence of irreducible
$SL_2({\mathbb C})$
representations of the fundamental group of a compact Riemannian manifold
is an equivariant harmonic
map to an ${\mathbb R}$-tree which is minimal and whose length function is
projectively
equivalent to the Morgan-Shalen limit of the sequence of representations.
We then examine the implications of the existence of a harmonic map when the
action on the tree fixes an end.
\end{abstract}

\maketitle

\section{Introduction}

For  a finitely generated group $\Gamma$, Morgan and Shalen \cite{MS}
compactified the character variety of equivalence classes of $SL_2({\mathbb
C})$
representations of
$\Gamma$ with  projective limits of the  length functions associated to the
representations.  These limits turn out to be projectively equivalent to
the length functions
of actions of
$\Gamma$ by isometries on ${\mathbb R}$-trees.  This built on earlier work of
Culler and Shalen \cite{CS} who had  identified the
ideal points of a complex curve of
$SL_2({\mathbb C})$ representations
 with  actions  on a simplicial tree.
The tree arose from  the Bass-Serre theory
applied to the function field of the curve with the discrete valuation
corresponding to the
point in question.

For unbounded sequences of discrete and faithful representations,
Bestvina \cite{B} and Paulin \cite{P} obtain an ${\mathbb R}$-tree whose
length function
is in the Morgan-Shalen class and which appears as the Gromov limit of
convex hulls
in ${\mathbb H}^3$. The limit involves rescaling the metric on the  hulls
by the maximum distortion at the center of the representation. Cooper
\cite{C} extends
Bestvina's  construction to obtain such a tree for sequences of representations
that are not necessarily discrete and faithful.  Moreover, he uses
length functions to show that if the sequence
eventually lies on a complex curve in the representation variety,
then the limiting tree is in fact simplicial.

In this paper we produce an ${\mathbb R}$-tree for any unbounded sequence
of irreducible
representations of the fundamental group of compact Riemannian manifolds along
with an equivariant harmonic map from the universal cover to the tree.
The starting point for this is to regard  representations of the
fundamental group as flat
connections on
$SL_2({\mathbb C})$  bundles, or equivalently, as harmonic maps from the
universal cover to ${\mathbb H}^3$. We first observe that
the non-existence of any convergent subsequence of a given sequence of
representations
 is equivalent to the statement that the
energy of the corresponding harmonic maps along all such subsequences is
unbounded. Then,
after rescaling by the energy, the recent work of Korevaar and Schoen
\cite{KS2}
applies: The harmonic maps pull back the metric from
${\mathbb H}^3$ to a sequence of (pseudo) metrics  which,
under suitable conditions, have a subsequence  that converges pointwise
to the pull back of a metric
on some non-positive curvature ($NPC$) space.
In our case, the conditions are met thanks to the Lipschitz property
of harmonic maps. We then show that the $NPC$ space is
an ${\mathbb R}$-tree with length function in the projective class of the
Morgan-Shalen
limit. In this way, rescaling by energy turns out to be strong enough to
give convergence,
but at the same time is subtle enough to give a non-trivial limiting length
function. The main result, which can also be thought of as an existence theorem
for harmonic maps to certain trees, is:

\medskip
\noindent
{\it The Korevaar-Schoen limit of an unbounded sequence of irreducible
$SL_2({\mathbb C})$
representations of the fundamental group of a compact Riemannian manifold
is an equivariant  harmonic map to an ${\mathbb R}$-tree. The image of this map
is a minimal subtree for the group action, and the projective class of
the associated length function  is the Morgan-Shalen
limit of the sequence.\it}
\medskip

\noindent
For harmonic maps into hyperbolic manifolds, Corlette \cite{C1}, \cite{C2}
and Labourie \cite{L} have shown that the
existence of a harmonic map implies that the action is semi-simple.
A partial analogue of this is our Theorem \ref{s-s} in the last section.
This also suggests that there is no obvious generalization for trees of
Hartman's result \cite{H} about uniqueness of harmonic maps.
Further results in this direction for surface groups are discussed in
\cite{DDW}.

We also note that this paper contains the work of Wolf \cite{W} for surfaces.
In the case considered there, sequences  of discrete, faithful $SL_2({\mathbb
R})$-representations give  equivariant harmonic maps ${\mathbb H}^2\to
{\mathbb H}^2$.  The limiting tree then arises as the leaf space
of the measured  foliation coming from the sequence of Hopf differentials.
A generalization of this point of view to $SL_2({\mathbb C})$  also appears
in \cite{DDW}.

\section{Background}

\noindent
{\bf A. Trees and Lengths:}
Recall that an ${\mathbb R}$-tree is a metric space
where any two points can be joined
by a unique arc, and this arc is isometric to an interval in ${\mathbb R}$.
For example, an increasing union of simplicial trees is an ${\mathbb R}$-tree.
Given a representation $\rho$ of a group $\Gamma$ to the isometry group of
an arbitrary
metric space
$X$, the length function is defined as
$$
  l_{\rho}(g)
  =
  \inf_{x \in X}
  d_X
  (x, \rho(g) x)  .
$$
In the case when $X$ is an ${\mathbb R}$-tree, we will use the fact  that
the length function
is  identically zero if and only $\Gamma$ has a fixed point (see \cite{CM}).
Also recall that length functions determine actions on trees, except
in some degenerate cases:
Suppose that two different actions of $\Gamma$ on two trees have the same,
non-abelian length function (\emph{non-abelian} means that $l(g)$ is not of
the form
$|h(g)|$ for some homomorphism
$h: \Gamma \to {\mathbb R}$). Then there is an equivariant isometry
between the respective (unique) minimal
subtrees of the same length function \cite{CM}.

\bigskip
\noindent
{\bf B. Character Variety:} Now take $\Gamma$ to be a finitely generated group.
 $\chi(\Gamma)$
will be its character variety, i.e. the space of characters of representations
of $\Gamma$ into $SL_2({\mathbb C})$.
Whereas the space of representations is merely an affine algebraic set
and $\chi(\Gamma)$ is a closed algebraic set,
the components of $\chi(\Gamma)$ containing
classes of irreducible representations are affine varieties
(closed, irreducible algebraic sets).
Conjugacy classes of irreducible representations
are characterized by their character (see \cite{CS}).

Given a representation $\rho: \Gamma \to SL_2({\mathbb C})$, for each
$g$ in $\Gamma$ the image $\rho(g)$ acts on ${\mathbb H}^3$ by
isometries. Consider the corresponding length function
$$
  l_{\rho}(g)
  =
  \inf_{x \in {\mathbb H}^3}
  d_{{\mathbb H}^3} (x, \rho( g ) x).
$$
 We may think of $l_{\rho}$ as a function on the generators
$\gamma_1,...,\gamma_r$
of $\Gamma$. 
For $C$ the set of conjugacy classes of $\Gamma$, the Morgan-Shalen
compactification of the character variety is
obtained by adding the projective limits of all
$\{ l_{\rho} (\gamma), \gamma \in C \}$
in the projective space
$\left( \prod [0, \infty) \setminus 0 \right) / {\mathbb R}^+$  (see
\cite{MS}).
Strictly speaking, \cite{MS} uses the traces to define the compactification.
That this is equivalent to using length functions follows from the fact that
if tr$\rho(g) \geq 1$ then
$|l_{\rho}(g) - 2 \ln$tr$\rho(g)| \leq 2$, see \cite{C}.
Explicitly, given a sequence of representations $\rho_n: \Gamma \to
SL_2({\mathbb C})$,
only one of the following can occur:
\begin{enumerate}
\item
For some subsequence $n'$, all traces $\rho_{n'}(\gamma_i)$ are bounded.
Then $\rho_{n'}$
converges (possibly after passing to a further
subsequence)  in $\chi(\Gamma)$ (see \cite{CL}).
\item
For every subsequence $n'$ there is some
$i$ such that tr$\rho_{n'}(\gamma_i) \to \infty$ as $n' \to \infty$. Then
there is an ${\mathbb R}$-tree and a representation $\rho: \Gamma \to Iso(T)$
such that $l_{\rho}$ is not identically zero and $l_{\rho_{n}} \to l_{\rho}$
projectively
(possibly after passing to a subsequence). 
\end{enumerate}

\bigskip
\noindent
{\bf C. Harmonic maps:} Recall from \cite{D} and \cite{C1} that
given a compact Riemannian manifold $M$ and an
irreducible representation $\rho$ of $\Gamma=\pi_1(M)$ into $SL_2({\mathbb
C})$,
any $\rho$-equivariant map from the universal cover $\widetilde M$ to
${\mathbb H}^3$
can be homotoped via the heat flow to a $\rho$-equivariant harmonic map
from the universal cover $\widetilde M$ to ${\mathbb H}^3$:
$$
  u: \widetilde M \to {\mathbb H}^3, \quad
  u(g x) = \rho(g) u(x), \quad
  D^* d u = 0,
$$
for $du:T \widetilde M \to T{\mathbb H}^3$ and $D$ the pull-back of the
Levi-Civita connection on ${\mathbb H}^3$.
This harmonic map minimizes the energy on $M$
$$
  \int_M
  | du |^2
  dM
$$
amongst all equivariant maps $v: \widetilde M \to {\mathbb H}^3$ (see p. 643 of
\cite{KS1})
and is unique up to
${\mathbb R}$-translations in some complete, totally geodesic submanifold
$Y \times {\mathbb R}$ of the target (see 4.4.B of \cite{GP}).

Conversely, given a $\rho$-equivariant harmonic map $u$ the representation
$\rho$ can be recovered as the holonomy of the flat connection
$A + i \Phi$,
for $A$ the pull-back by $u$ of the Levi-Civita connection on $T{\mathbb H}^3$
and $\Phi = -  du/2$. Hence,
$
E(u) =4 \|\Phi \|_2^2\
$, where $\Vert\cdot\Vert_2$ denotes the $L^2$ norm.
Flatness and harmonicity now become the equations:
\begin{eqnarray} \label{H}
  F_A = \frac{1}{2}[\Phi,\Phi] , \quad
  d_A \Phi = 0 , \quad
  d_A^* \Phi = 0.
\end{eqnarray}
We shall refer to $\Phi$ as the \emph{Higgs field} of the representation
$\rho$.
Thinking of representations as flat connections allows us to see easily that
a sequence escapes to infinity
only if the energy of $u$ blows up:
\begin{Prop}
Let $\rho_i$ be a sequence of representations of $\Gamma$ with Higgs fields
$\Phi_i$.
If
the energy  of the associated harmonic maps $u_i$ is bounded then there is
a representation
$\rho$  and a subsequence $i'$ such that
$\rho_{i'}
\to
\rho$ in  $\chi(\Gamma)$.
\end{Prop}
\begin{proof}
Recall that any harmonic map
$u: \widetilde M \to {\mathbb H}^3$ (or any target of non-positive curvature)
has the following property: For any $x$ in $\widetilde M$ and
$R>0$ there is
constant $C(R)$ independent of $u$ such that whenever $d(x,y) < R$,
\begin{eqnarray} \label{est}
  |\nabla u |^2 (y)
  \leq
  C(x,R)
  E_{B(x, R)}(u)  .
\end{eqnarray}
This follows from the Bochner formula for $|\nabla u|^2$
when the target has negative curvature, see \cite{S}.
Therefore $E(u_i) \leq B$ implies uniform bounds on $\| \Phi_i \| _{C_0} $.
By the first of
equations (\ref{H}) this implies uniform bounds on $\| F_{A_i} \|_{C_0}$
and hence
uniform bounds on $\| F_{A_i} \|_p$ for any $p$. 
By standard application of Uhlenbeck's weak compactness theorem
and elliptic regularity the result follows.
\end{proof}

\bigskip
\noindent
{\bf D. Korevaar-Schoen compactness:}
The previous section shows that in order to examine the ideal points of
the character variety one has to look at equivariant harmonic maps of
arbitrarily high energy. For this, we recall the construction in \cite{KS2}.

Let $\Omega$ be a set and $u$ a map into an $NPC$ space $X$.
Use $u$ to define the pull-back pseudometric on $\Omega$,
$
  d_u (x, y)
  =
  d_{X}
  ( u(x), u(y) )
$,
for any $x$ and $y$ in $\Omega$. To achieve convergence in an $NPC$ setting,
some convexity is needed.  Korevaar and Schoen achieve this by enlarging
$(\Omega, d_u)$ to a space $\Omega_{\infty}$ obtained by adding the segments
joining any two points in $\Omega$,
 the segments joining any two points on these segments,
and so on. Then they extend  the pull-back pseudometric  from $\Omega$
to $\Omega_{\infty}$ linearly.
After identifying
points of zero pseudodistance in $(\Omega_{\infty}, d_u)$ and completing,
one obtains a metric space $(Z,d_u)$ isometric
to the convex hull $C(u(\Omega))$ in the target $X$ (and hence $NPC$).

It is a crucial point that certain inequalities (which carry over to pointwise
limits) satisfied by pull-back pseudodistances are enough for this $Z$ to
be an $NPC$ space,
regardless of the pseudodistance being a pull-back or not
(see Lemma $3.1$ of \cite{KS2}).
The following
summarizes the main results from \cite{KS2} needed here when
$\Omega = \widetilde M$ is the universal cover of a compact Riemannian
manifold:

\begin{Thm} \label{KS2}
Let $u_k: \widetilde M \to X_k$ be a sequence of maps on the universal cover of
 some Riemannian manifold $\widetilde M$
such that:
\begin{enumerate}
\item
Each $X_k$ is an $NPC$ space;
\item
The $u_k$'s have uniform modulus of continuity: For each $x$ there is a
monotone
function $\omega(x,)$ so that
$\displaystyle
  \lim_{R \to 0} \omega(x, R) = 0\ $, and
$\displaystyle
  \max_{B(x,R)} d_{u_k}(x,y) \leq \omega(x, R)\
$.
\end{enumerate}
Then
\begin{enumerate}
\item \label{conv}
The pull-back
pseudometrics $d_{u_k}$ converge (possibly after passing to a subsequence)
pointwise,
locally uniformly, to a pseudometric
$d_{\infty}$;
\item \label{NPC}
The Korevaar-Schoen construction for $(\widetilde M, d_{\infty})$ produces
a metric space
$(Z, d_{\infty})$ which is $NPC$;
\item \label{min}
If the $u_k:\widetilde M \to X_k$ have uniformly bounded energies and are
energy minimizers
then the projection $u:(\widetilde M, d_{\infty}) \to Z$ is also energy
minimizer;
\item \label{equiv}
If the $u_k$'s are equivariant then $u$ is also equivariant;
\item \label{nontr}
If $\displaystyle \lim_{k \to \infty} E(u_k)$ is not zero then  $u$ is not
trivial.
\end{enumerate}
\end{Thm}
\begin{proof}
\ref{conv}. and \ref{NPC}. are contained in
Proposition $3.7$ of \cite{KS2}.
Given the locally uniform convergence, \ref{min}., \ref{equiv}. and
\ref{nontr}. follow from Theorem $3.9$ of
\cite{KS2}.
\end{proof}

\section{Main Results}

In this section, we let
 $\Gamma$ be the fundamental group of a compact Riemannian manifold $M$,
and let $\rho_k$ be a sequence of irreducible representations
of $\Gamma$ in $SL_2({\mathbb C})$ with no convergent subsequence
in $\chi(\Gamma)$ (henceforth, an \emph{unbounded} sequence).
Let $u_k: \widetilde M \to {\mathbb H}^3$ be the corresponding harmonic maps.
Rescale the metric $d_{{\mathbb H}^3}$ on  ${\mathbb H}^3$ to
$$
  {\hat d}_{{\mathbb H}^3, k}
  =
  \frac{d_{{\mathbb H}^3}}{E(u_k)^{1/2} }
$$
before pulling it back via $u_k$ to
$
  {\tilde d}_k (x, y)
  =
  {\hat d}_{{\mathbb H}^3,k} (u_k (x), u_k(y))
$
on $\widetilde M$.
Continue to denote the harmonic maps into the rescaled targets by $u_k$.
With this understood, the first result of this paper may be stated as
follows:

\begin{Thm}  \label{main}
Let $M$ be a compact Riemannian manifold
and $\rho_k$  an unbounded sequence of irreducible $SL_2({\mathbb C})$
representations of
$\Gamma=\pi_1(M)$.
Then the Korevaar-Schoen limit of the rescaled harmonic
maps
$
  u_k: \widetilde M \to \left( {\mathbb H}^3, {\hat d}_{{\mathbb H}^3,k}
\right)
$
is an energy minimizing $u:\widetilde M \to T$,
for $T$ an ${\mathbb R}$-tree. In addition, $\Gamma$ acts on $T$ without
fixed points,
and $u$ is $\Gamma$-equivariant.
\end{Thm}
\begin{proof}
According to \ref{conv}. of Theorem \ref{KS2},
to show convergence of the $u_k$'s it is enough to show
that they have uniform modulus of continuity.
Recall the estimate (2).
With this, the rescaled sequence $u_k$
satisfies the uniform modulus of continuity condition for $\omega(x, R) =
R$ for all $x$:
$
  {\tilde d}_k (u_k(x), u_k (y))
  \leq
  d_M(x,y)\
$.
Therefore the
pseudodistances ${\tilde d}_k$  converge   pointwise and locally uniformly
to a limiting pseudodistance $\tilde d$ on $\widetilde M$.
Let $u$ be the projection
$
  u:\widetilde M \to Z\
$
for $Z$ the $NPC$ metric space obtained
from $(\widetilde M, \tilde d)$ by identifying points of zero distance and
completing.
Then $u$ is energy minimizing by \ref{min}. of Theorem \ref{KS2}
and the harmonicity of the $u_k$'s.

To show that  the limiting $NPC$ space $Z$ is in fact
an ${\mathbb R}$-tree, we need to show that:
a) Any two points in $Z$ can be joined by a unique arc;
b) Every arc in $Z$ is isometric to a segment in ${\mathbb R}$.
Part b) follows immediately from the defining property of $NPC$ spaces
(any two points
can be joined by an arc isometric to a segment in ${\mathbb R}$).
For a) we must show that there are no more arcs.
If $d_k$ is the pull-back of the standard metric on ${\mathbb H}^3$
then for each $k$ the metric space
$
  Z_k = \left( \widetilde M _{\infty}, d_{k} \right) / \sim\
$
is, by construction, isometric to the convex hull of the image
$u_k(\widetilde M)$ in ${\mathbb
H}^3$.  Therefore, for any two points $x,y\in\widetilde M _{\infty}$ the
geodesic
segment isometric to $[u_k(x),u_k(y)]$ lies in $Z_k$. Suppose $z\neq x,y$
is a third point
on  some other arc joining
$x$ and $y$, and consider the geodesic triangle
with vertices $u(x)$, $u(y)$, and $u(z)$ in ${\mathbb H}^3$.
It is a standard fact that there is a constant $C$ characteristic of
${\mathbb H}^3$
such that in the standard metric $d_{{\mathbb H}^3}$
any point in the interior of this triangle has distance less than $C$ from the
edges (see the proof of Theorem $3.3$ of \cite{B}).
Then in the rescaled pull-back metric ${\tilde d}_k$,
any point in the  interior has distance less than $C/E(u_k)^{1/2}$; hence,
as $E(u_k) \to \infty$ this distance becomes
arbitrarily small.
Therefore, all triangles become infinitely thin at the limit.
This suffices to show that there can be only one arc joining any two points,
see Proposition  $6.3.$C of \cite{Gr} and page $31$ of \cite{GH}.

In addition, because of the rescaling $E(u_k) = 1$ for all $k$ and using
\ref{nontr} of
Theorem \ref{KS2},
 the limit $u$ is also non-trivial.
Therefore, $T=Z$ is non-trivial.

The action of $\Gamma$ on $\widetilde M$ extends to an action on the whole of
$\widetilde M _{\infty}$ (this follows from a straightforward calculation
using the fact that
on each segment $[x,y]$ of length $d$ there is only one point
 $\lambda d$ away  from $x$ and $(1 -\lambda)d$ away from $y$).
The equivariance of each $u_k$ and $u$ implies that there are actions
 $\sigma_k$ on $( \widetilde M_{\infty},\tilde d_k )$ and $\sigma$ on
$(\widetilde M _{\infty},\tilde d)$
 by isometries.
These clearly descend to the completed quotients $Z_k$ and $T$.
Now suppose that $\Gamma$ acted on $T$ with some fixed point $t_0$ .
Then the constant map $w(x) = t_0$
is equivariant with respect to $\sigma$ and of zero energy.
Then $u$, the energy minimizer, also has zero energy; a contradiction.
Hence, there are no  points on $T$ fixed by all elements of
$\Gamma$, which is equivalent to the length function $l_s$ being non-zero.
\end{proof}

The pointwise convergence of the pseudometrics as $k \to \infty$
implies that the length function
$$
  l_{\sigma}(g)
   =
  \inf_{x \in \Omega_{\infty}}
  d_{{\mathbb H}^3} (  u(x), u(gx)  )
$$
is the pointwise limit of the length functions
$$
  l_{\sigma _k}(g)
  =
  \inf_{x \in \Omega_{\infty}}
  d_{{\mathbb H}^3} (  u_k(x), u_k(gx)  ).
$$
On the other hand the $u_k$'s came from representations $\rho_k$, so let
$l_{\rho_k}$
be the length function obtained by the action of the $\rho(g)$'s on
${\mathbb H}^3$
$$
  l_{\rho_k}(g)
  =
  \inf_{v \in {\mathbb H}^3}
  d_{{\mathbb H}^3}  ( v, \rho_k(g) v ) ,
$$
Recall that the Morgan-Shalen limit of the $\rho_k$'s in $\chi(\Gamma)$
is the projective limit of the length functions $l_{\rho_k}$.

\begin{Thm} \label{lengths}
The length function $l_{\sigma}$ of the action of $\Gamma$ on the
Korevaar-Schoen
 tree $T$ is in the projective class
of the Morgan-Shalen limit of the sequence $\rho_k$.
\end{Thm}

\begin{Cor}
For a length function $l$ appearing as the Morgan-Shalen limit of
irreducible elements of $\chi(\Gamma)$, there is an ${\mathbb R}$-tree $T$
on which $\Gamma$ acts by isometries with
length function $l$, and  an equivariant harmonic map $u:\widetilde M \to T$.
\end{Cor}

\begin{proof}[Proof of Theorem \ref{lengths}]
Since
$
  l_{\sigma}
  =
  \lim_{k \to \infty}
  l_{\sigma_k}\
$,
we only need to show that $l_{\rho_k}$ and $l_{\sigma_k}$ converge
projectively to the same (non-trivial) limit.
But
$$
  l_{\sigma _k}(g)
  =
  \inf_{x \in \Omega_{\infty}}
  d_{{\mathbb H}^3} (  u_k(x), u_k(gx)  )
$$
where the $\inf$ on the right hand side is over
the lengths of the geodesics in ${\mathbb H}^3$
joining $u_k(x)$ to  $u_k(gx)$.
Now according to Lemma $2.5$ of \cite{C}, and as a result
of the property of thin triangles in ${\mathbb H}^3$,
this  geodesic  contains a subgeodesic
with end points $A$ and $B$ such that:
$$
  \left| \left|[A,B]\right| - l_{\rho_k}(g) \right|
  \leq
  \Delta\ ,\qquad
  d_{{\mathbb H}^3} ( B, \rho_k(g) A)
  \leq
  \Delta\ .
$$
This implies
$
  d_{{\mathbb H}^3} ( A, \rho_k(g) A)
  \leq
  l_{\rho_k}(g) + 2 \Delta\
$.
Now by construction there is $x\in\Omega_{\infty}$ such that $u_k (x) = A$, and
therefore
$
  l_{\sigma_k} (g)
  \leq
  l_{\rho_k}(g) + 2 \Delta\
$.
It also follows from the definitions that $l_{\rho_k} (g) \leq l_{\sigma_k}
(g)$.
Dividing each side by $E(u_k) \to \infty$, we get the same limit.
Since the action of $\Gamma$ on $T$ has no fixed points, this limit is
non-trivial.
This completes the proof of Theorem \ref{lengths}.
\end{proof}

\section{Minimality of the tree}

The purpose of this section is to show that the image of the equivariant
harmonic map
of the previous section is a minimal tree,
i.e. it does not contain any proper subtree invariant under the action of
$\Gamma$. The main idea is that there is not enough energy for such a
subtree.

To begin, recall that given a closed subtree $T_1$ of an $\mathbb R$-tree
$T$ and a
point $p$  not in $T_1$, there is a unique shortest arc from $p$ to $T_1$,
obtained as
the closure of $\gamma \setminus T_1$ for any arc $\gamma$ from $p$ to
$T_1$ (see
$1.1$ of \cite{CM}). Call
the endpoint of this unique arc $\pi(p)$.

\begin{Lem}
The map $\pi$ is distance decreasing.
\end{Lem}

\begin{proof}
Given two points $p$ and $q$ not on $T_1$, the unique arc in $T$
joining them
either contains $\pi(p)$ or not.
If it does, then it also contains $\pi(q)$, since $T_1$ is a subtree.
Therefore $ d( \pi(p), \pi(q) ) < d( p,q )$.
If it doesn't, then $d( \pi(p), \pi(q) ) = 0$.
\end{proof}

\begin{Lem}  Suppose $\Gamma$ acts on $T$.
If $T_1$ is $\Gamma$-invariant, then
$\pi$ is $\Gamma$-equivariant.
\end{Lem}
\begin{proof}
To prove equivariance note that by invariance of the subtree, if $p$
is not in $T_1$ then $gp$ is not,
either. Then if $\pi(gp) \neq g\pi(p)$, by definition
$
  d(\pi(gp),gp)
  <
  d(g\pi(p),gp)
$,
and since $g$ is an isometry
$
  d(g^{-1} \pi(gp),p)
  <
  d(\pi(p),p)
$,
which contradicts the definition of $\pi(p)$.
\end{proof}

Now suppose that $u:\widetilde M\to T$ is an equivariant harmonic map.

\begin{Lem}   \label{subtree}
If $T_1$ is a proper subtree of $T$ invariant under the action of $\Gamma$,
and $u(\widetilde M) \cap T_1 \neq \emptyset$,
then $u(\widetilde M) \subset T_1$.
\end{Lem}
\begin{proof}
Suppose that there is a point $u(x)$ not in $T_1$. The image of a closed ball
around $x$
in $\widetilde M$, large enough to enclose a fundamental domain,
consists of finitely many arcs starting from $u(x)$ in $T$.
Therefore, there is $\varepsilon$
such that for any $w_1$ and $w_2$ in the image of $u$ and
within distance $\varepsilon$ of $u(x)$
we have
$$
  d(\pi(w_1), \pi(w_2))
  \leq
  \frac{1}{2}
  d(w_1,w_2).
$$
Then by the definition of the energy density
$$
  |du|^2 (x)d\mu
  =
  \lim_{\varepsilon \to 0}
          \frac{1}{\omega_n}
          \int_{S(x,\varepsilon)}
               \frac{d^2(u(x),u(y))}{\varepsilon^2}
          \frac{d\sigma(y)}{\varepsilon^{n-1}}
$$
(see \cite{KS2}, pp. 227-228)
we have $|d( \pi \circ u)| d \mu  \leq \frac{1}{2} |d u| d \mu  $
on $u^{-1}\left(  B_{u(x)} (\varepsilon) \right)$.
Since $u$ is an energy minimizer, this implies that $u$ is constant.
Therefore, the image of $u$ cannot intersect $T_1$; a contradiction.
\end{proof}

As a consequence, we have the following:

\begin{Thm}
Let $M$ be a compact Riemannian manifold. 
For $u:\widetilde M\to T$ a $\Gamma$-equivariant harmonic map,
the image of $u$ is a minimal subtree of $T$.
In particular, the image of $u$ in the Korevaar-Schoen tree in
Theorem \ref{main} is minimal.
\end{Thm}

\begin{Cor}
Suppose that the sequence of representations lies on a complex
curve in the character
variety. Then its Korevaar-Schoen limit is a simplicial tree.
\end{Cor}
\begin{proof}
If the sequence or representations lies on a curve then the limit tree $T$
contains an invariant simplicial subtree $T_0$, see \cite{CS}, \cite{C}.
The corollary follows by the minimality of the image of $u$ and the uniqueness
of minimal subtrees (Proposition $3.1$ of \cite{CM}).
\end{proof}

\section{Actions with fixed ends}

Recall that a  {\it ray} in a tree $T$ is an arc isometric to $[0, +
\infty)$ in
${\mathbb R}$.
Two rays are said to be equivalent if their intersection is still a ray,
and an equivalence
class of rays is an {\it end}. An end is  fixed under $\Gamma$ if
$g R \cap R$ is a ray for any
$g$ in $\Gamma$ and ray $R$ in the end.  An action of $\Gamma$ on an
$\mathbb R$-tree $T$ is said to be \emph{semi-simple} if either $T$ is
equivariantly isometric to an action on $\mathbb R$, or $\Gamma$ has
no fixed ends.

Suppose that the action of $\Gamma$ is such that it fixes some end.
Given a ray $R$ in the fixed end,
let $\pi(x)$ denote the projection of any point $x$ to $R$
as above, and let $R_{\pi(x)}$ be the sub-ray of $R$ starting from
$\pi(x)$.
Then $R_x = [x, \pi(x)] \cup R_{\pi(x)}$ is the unique ray starting
from any given point
$x$ and belonging to the given end.
Now define $\phi_\varepsilon: T \to T$
by
taking $\phi_\varepsilon(x)$ to be the unique point of distance
$\varepsilon$ from $x$
on $R_x$.

\begin{Lem}
$\phi_\varepsilon$ is equivariant.
\end{Lem}
\begin{proof}
Since the end is
invariant, $g R_x = R_{gx}$, and $R_{gx}$ starts from $gx$. Now the point
$g \phi_\varepsilon(x)$ is on $R_{gx}$ and is a distance $\varepsilon$ from
$gx$.
Then by definition $\phi_\varepsilon(gx) =
g \phi_\varepsilon(x)$.
\end{proof}

\begin{Lem} \label{dd}
$\phi_\varepsilon$ is distance decreasing.
\end{Lem}
\begin{proof}
For $x$ and $y$  with $R_x$ subray of $R_y$, $\phi_\varepsilon$ is shift by
distance $\varepsilon$
along the ray and $d(\phi_\varepsilon(x), \phi_\varepsilon(y)) = d(x,y)$.
If $R_x$ and $R_y$ intersect along a proper subray of both, let $p$ be the
initial point of this subray. Then the arc from $x$ to $y$
contains
$p$ and $d(x,y) = d(x,p) + d(p,y)$.
There are three cases to consider:
\begin{enumerate}
\item
$d(x,p) \geq \varepsilon$, $d(y,p) \geq \varepsilon$. Then
$d(\phi_\varepsilon(x),\phi_\varepsilon(y)) = d(x,p) -\varepsilon +
d(y,p) - \varepsilon  \leq d(x,y)$.
\item
$d(x,p) \leq \varepsilon$, $d(y,p) \leq \varepsilon$.  Then
$d(\phi_\varepsilon(x),\phi_\varepsilon(y)) = |d(x,p) -\varepsilon - (
d(y,p) - \varepsilon ) |  \leq  d(x,y)$.
\item
$d(x,p) \leq \varepsilon \leq d(y,p)$ (or symmetrically,
$d(y,p) \leq \varepsilon \leq d(x,p)$).
Then $d(\phi_\varepsilon(x),\phi_\varepsilon(y)) =
(\varepsilon - d(x,p) )  + ( d (y,p) - \varepsilon) = d(y,p) - d(x,p)
\leq d(x,y)$.
\end{enumerate}
\end{proof}

\noindent
Now, arguing as in the proof of Lemma \ref{subtree} applied to
$\phi_\varepsilon \circ
u$, and using Lemma
 \ref{dd}, we conclude the following:

\begin{Thm} \label{s-s}
Let $M$ be a compact Riemannian manifold, $T$ an ${\mathbb R}$-tree on which
$\Gamma=\pi_1(M)$ acts minimally and non-trivially via isometries, and suppose
that there is a  $\Gamma$-equivariant energy minimizing map $u:
\widetilde M
\to T$. Then either the action of $\Gamma$ on $T$ is semi-simple, or $u$ is
contained in
a continuous family of distinct $\Gamma$-equivariant energy minimizers.
\end{Thm}

\end{document}